\documentclass[reqno]{amsart}
\usepackage{amssymb,amsthm,amsfonts,amstext}
\usepackage{amsmath}

\makeatletter \@addtoreset{equation}{section} \makeatother

\renewcommand\thefigure{\thesection.\@arabic\c@figure}
\renewcommand\thetable{\thesection.\@arabic\c@table}
\newtheorem{theorem}{Theorem}[section]
\newtheorem{lemma}[theorem]{Lemma}
\newtheorem{proposition}[theorem]{Proposition}

\newtheorem{remark}[theorem]{Remark}

\newcommand{\mc}[1]{{\mathcal #1}}

\newcommand{\bb}[1]{{\mathbb #1}}

\newcommand{\<}{\langle}
\renewcommand{\>}{\rangle}

\author{Milton D. Jara V.}
\date{}
\title[Invariance principle for a tagged particle]{Non-equilibrium scaling limit
for a tagged particle in the simple exclusion process with long jumps}

\begin{document}

\begin{abstract}
We prove an invariance principle for a tagged particle in a simple exclusion
process with long jumps out of equilibrium. 
\end{abstract}


\renewcommand{\subjclassname}{\textup{2000} Mathematics Subject Classification}

\address{IMPA, Estrada Dona Castorina 110, CEP 22460-320, Rio de Janeiro,
Brazil} 
\email{monets@impa.br}

\maketitle
\section{Introduction}

A classical problem in statistical mechanics consists in proving that the
dynamics of a single (tagged) particle in a particle system satisfies an
invariance principle.  
The first important result on the position of a tagged particle in the diffusive
scaling is due to Kipnis and Varadhan \cite{KV}. In a remarkable paper, they
develop an invariance principle for additive functionals of a reversible Markov
process and deduce a central limit theorem for a tagged particle in symmetric
simple exclusion processes. The approach requires the process to be in
equilibrium, that is, that the initial distribution of particles corresponds to
an equilibrium measure. This result is robust and can be extended to other
situations \cite{Var}, \cite{SVY} (see also \cite{Fer}, \cite{LOV} and the
references therein). 

Out of equilibrium, the asymptotic behavior of a tagged particle appears as one
of the central problems in the theory of interacting particle systems, and
remains mostly unsolved. Even a law of large numbers for a tagged particle on a
system starting from Bernoulli product measures of slowly varying parameter
seems out of reach. Notice however the result of Rezakhanlou \cite{Rez} stating
that the average behavior of tagged particles is given by a diffusion process.  
Some results in dimension $d=1$ have been obtained \cite{JL1}, \cite{JL2}. These
results are based on the following observation. For nearest-neighbor systems in
$d=1$, the initial ordering of the particles is preserved, and the position of a
tagged particle can be related to the empirical density of particles and the
total current of particles through the origin. 

Recently, a general approach to this problem has been proposed \cite{JLS}. The
idea is to obtain the hydrodynamic limit of the process {\em seen by an observer
sitting on the tagged particle}. Let us call this process the {\em environment}
process. The density of particles for the environment process can be obtained by
shifting the density of particles of the original process to the position of the
tagged particle. A martingale representation of the position of the tagged
particle should allow to obtain the scaling limit of the tagged particle as the
solution of some martingale problem. This strategy has been completed in
\cite{JLS} for a mean-zero, zero-range process in dimension $d=1$.  

However, the position of a tagged particle is a martingale for the zero-range
process. We already know from \cite{Rez} that the scaling limit of a tagged
particle is not in general a martingale. Therefore, it is not clear if this
strategy can be followed for the simple exclusion process.  

In this article, we consider an exclusion process with long-range jumps. This is
a simple exclusion process with transition rate $p(x,y) = |y-x|^{-(d+\alpha)}$,
$\alpha \in (0,2)$. We obtain a non-equilibrium invariance principle for the
position of a tagged particle. The limiting process is a time-inhomogeneous
process of independent increments, driven by a solution $u(t,x)$ of the
fractional heat equation $\partial_t u =\Delta^{\alpha/2}$, where
$\Delta^{\alpha/2}$ is the fractional Laplacian. In the equilibrium case, the
limiting process is a symmetric, $\alpha$-stable L\'evy process.  

In \cite{KL}, the position of the tagged particle is written as a martingale
plus a an additive functional of the environment process. The invariance
principle is therefore obtained as a consequence of the central limit theorem
for martingales and for additive functionals of reversible Markov processes. 
This approach can not be followed for the exclusion process with long jumps,
since the tagged particle does not have bounded moments of order greater than
$\alpha$, and in particular, of order two. We therefore represent complex
exponentials of the tagged particle as the product of exponential martingales
and additive functionals of the environment process. It turns out that, due to
the long range interaction, only a law of large numbers for additive functionals
of the environment process is needed. This law of large numbers is a consequence
of the hydrodynamic limit for the environment process, and the strategy of
\cite{JLS} can be carried out. 
 
The main feature that allow us to solve the tagged particle problem for this
model, is that the environment process satisfies the {\em gradient condition}
(\cite{KL}, Chapter 5). This condition is not satisfied by the symmetric simple
exclusion process of finite range considered in \cite{KV}. This gradient
condition is equivalent to the fact that only a law of large numbers for
additive functionals of the environment process is needed to obtain the scaling
limit of the tagged particle. Notice that the non-gradient method \cite{VY} has
been developed to overcome this difficulty for the usual hydrodynamic limit of
systems that does not satisfy the gradient condition. A natural open problem is
to extend the non-gradient method to this setting, to obtain an invariance
principle for  the tagged particle problem out of equilibrium for the simple
exclusion process with finite-range jumps. Notice that when $1 \leq \alpha <2$,
we already need some non-gradient tools, to get rid of a non-gradient term that
does not contribute to the limit.

For the particular case of the simple exclusion process with long jumps, the
additive functional appearing in the martingale decomposition of the tagged
particle is already a (singular) function of the empirical measure. For general
processes with local interaction, a local replacement lemma confines us to
dimension $d <2$.   

The structure of the article is the following. In Section \ref{s4} we give
precise definitions of the simple exclusion process with long jumps and the
environment process, we state the main results and we define the hydrodynamic
limit of the process.  
In Section \ref{s5} we obtain the scaling limit of the tagged particle in
equilibrium. Although this result is a consequence of the non-equilibrium result
of Section \ref{s6}, we first give a complete proof of the equilibrium result
for two reasons. We obtain in this way a simplified exposition, and the proof of
the equilibrium invariance principle can be easily generalized to particle
systems with interactions and long jumps. Notice that for the case of a
zero-range process with long jumps, the fact that the tagged particle is already
a martingale does not really simplify the proof in the equilibrium case, since
we are not aware of a characterization of a L\'evy process avoiding exponential
martingales, that could be of use in this situation. 
In Section \ref{s6}, we obtain the scaling limit for the tagged particle out of
equilibrium, characterizing the limiting points by means of a martingale
problem.

\section{The exclusion process with long jumps}
\label{s4}
The simple exclusion process with long jumps is a system of long-range random
walks 
on the lattice $\bb Z^d$, conditioned to have at most one particle per site.
Consider $p: \bb Z^d \to \bb R_+$, the transition probability of a simple random
walk in $\bb Z^d$. The dynamics of this process can be described as follows. A
particle at site $x \in \bb Z^d$ waits an exponential time of rate 1, at the end
of which the particle tries to jump to site $x+z$ with probability $p(z)$. If
the site $x+y$ is empty, the jump is accomplished. Otherwise, the particle stays
at $x$ and a new exponential time starts. This is done independently for each
particle.  

This dynamics corresponds to a Markov process $\eta_t$ defined on the state
space $\Omega = \{0,1\}^{\bb Z^d}$. For a configuration $\eta \in \Omega$,
$\eta(x)=1$ represents a particle at site $x \in \bb Z^d$ and $\eta(x)=0$
represents an empty site. The process $\eta_t$ is generated by the operator  
\[
L f(\eta) = \sum_{x, y \in \bb Z^d} p(y-x)
\eta(x)\big(1-\eta(y)\big)\big[f(\eta^{xy})-f(\eta)\big], 
\]
where $f: \Omega \to \bb R$ is a function that depends on a finite number of
coordinates, $\eta$ is an element of $\Omega$ and $\eta^{xy}$ is the
configuration obtained from $\eta$ by exchanging $\eta(x)$ and $\eta(y)$: 
\[
\eta^{xy}(z)=
\begin{cases}
\eta(y), z=x \\
\eta(x), z=y\\
\eta(z), z \neq x,y.
\end{cases}
\]

The set of local functions in $\Omega$ is a core for the operator $L$
\cite{Lig}. We say 
that the transition rate $p(\cdot)$ is homogeneous, regular of degree $\alpha$
if there exists a function $q: \bb R^d \setminus \{0\}$ of class $\mc C^2$ such
that $p(z) = q(z)$ for any $z \in \bb Z^d \setminus \{0\} \to \bb R$ and such
that 
$q(\lambda u) = \lambda^{d+\alpha}q(u)$ for any $\lambda \neq 0$ and any $u \in
\bb R^d \setminus \{0\}$. Without loss of generality, we assume  $p(0)=0$.
The canonical example of such a rate is $p(z) =
p^*|z|^{-(d+\alpha)}$, where $p^*$ is the normalizing constant and $|\cdot|$ is
the Euclidean norm. Since $\sum_z p(z) < +\infty$, we have the restriction
$\alpha >0$. Since we are interested in non-diffusive scaling limits, we also
assume $\alpha <2$. The asumption $\lambda \neq 0$ instead of $\lambda >0$
restricts ourselves to {\em symmetric} homogeneous functions $q$.  

For each $\rho \in [0,1]$, denote by $\mu_\rho$ the product measure in $\Omega$
such that  
\[
 \mu_\rho(\eta(x)=1)=1-\mu_\rho(\eta(x)=0)= \rho.
\]
Due to the translation invariance of the transition rates, the measures
$\{\mu_\rho; \rho \in [,1]\}$ are invariant and ergodic for the process
$\eta_t$.

\subsection{The tagged particle}

Let $\eta \in \Omega$ be an initial configuration of particles such that
$\eta(0)=1$. We follow the evolution of the particle initially at the origin
together with the evolution of $\eta_t$. We call this particle the {\em tagged
particle} and we denote by $X_t$ its position at time $t$. We have the following
invariance principle for $X_t$: 

\begin{theorem}
\label{t1}
 Let $\eta_t$ be the simple exclusion process corresponding to an homogeneous,
 regular transition rate of degree $\alpha$. Fix a density $\rho \in [0,1]$ and
 assume that the process $\eta_t$ starts from the measure
 $\mu_\rho(\cdot|\eta(0)=1)$. Then the process $X_t^n=: n^{-1} X_{tn^\alpha}$
 converges in distribution on $\mc D([0,\infty],\bb R^d)$ to the L\'evy process
 $(1-\rho)Z_t$, where $Z_t$ is characterized by 
\[
 -\log E[\exp\{i\beta Z_t\}] = t\psi(\beta), \beta \in \bb R^d,
\]
 \[
  \psi(\beta) = \int_{\bb R^d} \big(1-e^{i\beta u}\big)q(u) du,
 \]
where we take the principal Cauchy value of the integral in the definition of
$\psi(\beta)$. 
\end{theorem}

In order to prove this result, it will be convenient to introduce an auxiliary
process (see \cite{KV}). We define $\xi_t$ by the relation $\xi_t(z)=\eta_t(X_t
+z)$, $z \in \bb Z^d$. With this definition, $\xi_t(0) \equiv 1$, so we consider
$\xi_t$ as a process in $\Omega_*= \{0,1\}^{\bb Z^d_*}$, where $\bb Z^d_* = \bb
Z^d \setminus \{0\}$. We call $\xi_t$ the {\em environment as seen by the tagged
particle}. Notice that both $\xi_t$ and $(\eta_t,X_t)$ are Markov processes. The
advantage in considering $\xi_t$ instead of $(\eta_t,X_t)$ is that the former
admits invariant measures, while the latter does not. In fact, for any $\rho \in
[0,1]$, the product measure $\nu_\rho$ in $\Omega_*$ given by 
\[
 \nu_\rho(\eta(x)=1)= 1-\nu_\rho(\eta(x)=0) = \rho \text{ for any } x \in \bb
 Z^d_* 
\]
is ergodic and invariant for $\xi_t$. The generator of $\xi_t$ is given by 
\[
 \mc L f(\xi) = \sum_{x, y \in \bb Z^d_*} p(y-x) [f(\xi^{xy})-f(\xi)]
	+ \sum_{z \in \bb Z^d_*} p(z)\big(1-\eta(z)\big)[f(\theta_z
	\xi)-f(\xi)], 
\]
where $\theta_z \xi$ corresponds to the configuration
\[
 \theta_z \xi(z) =\begin{cases}
                   \xi(x+z), &x \neq -z,0\\ \xi(z), &x = -z.
                  \end{cases}
\]

The position $X_t$ of the tagged particle can be recovered from the evolution of
$\xi_t$ as follows. Let $N_t^z$ be the number of translation of $\xi_t$ in
direction $z$ up to time $t$. Clearly, 
\[
 X_t = \sum_{z \in \bb Z^d_*} zN_t^z.
\]

\subsection{The hydrodynamic limit}

Let $u_0: \bb R^d \to [0,1]$ be a continuous function. We say that a 
sequence $\{\mu^n\}_n$ of probability measures in $\Omega$ is associated to
$u_0$  
if for any $\epsilon >0$ and any continuous function $G: \bb R^d \to \bb R$ of 
bounded support we have
\[
\lim_{n \to \infty} \mu^n\big\{\eta \in \Omega; 
\big|n^{-d} \sum_{z \in \bb Z^d} \eta(z) G(z/n) - \int G(x) u_0(x) dx\big| >
\epsilon\big\}=0. 
\]

Fix a sequence $\{\mu^n\}_n$ of measures in $\Omega$.   
We denote by $\eta_t^n$ the process $\eta_{tn^\alpha}$ starting from $\mu^n$.
The distribution of $\eta^n_t$ in $\mc D([0,\infty),\Omega)$ is denoted by $\bb
P^n$ and the expectation with respect to $\bb P^n$ is denoted by $\bb E^n$.  
We define the empirical density of particles $\pi_t^n$ by
\[
\pi_t^n(dx) = \frac{1}{n^d} \sum_{z \in \bb Z^d} \eta_t^n(z) \delta_{z/n}(dx),
\]
where $\delta_x$ is the Dirac-$\delta$ distribution at $x \in \bb R^d$. Notice
that $\pi_t^n$ is a process in $\mc D([0,\infty), \mc M_+(\bb R^d))$, where $\mc
M_+(\bb R^d)$ is the set of positive, Radon measures in $\bb R^d$. 
The following result can be proved as in (\cite{KL}, Chapter 4):

\begin{theorem}
\label{t3}
Let $u_0$ be a continuous, integrable initial density profile. Assume that
$\{\mu^n \}_n$ is associated to $u_0$.  
Then, the process $\pi_t^n$ converges in probability to the deterministic path
$u(t,x) dx$,  
where $u(t,x)$ is the solution of the hydrodynamic equation
\begin{equation}
\label{echid}
\begin{cases}
\partial_t u &= \bb L u \\
u(0,\cdot) &= u_0(\cdot)
\end{cases}
\end{equation}
and $\bb L$ is the integral operator given by
\[
\bb L F(x) = \int_{\bb R^d} q(y)\big\{F(x+y)+F(x-y)-2F(x)\big\}dy.
\]

In particular, for any time $t>0$, any $\epsilon >0$ and any continuous function
$G: \bb R^d \to \bb R$ of bounded support we have 
\[
\lim_{n \to \infty} \bb P^n\big\{
	\big|n^{-d} \sum_{z \in \bb Z^d} \eta_t^n(z) G(z/n) - \int G(x) u(t,x)
	dx\big| > \epsilon\big\}=0, 
\]
\end{theorem}

In order to obtain an invariance principle for $X_t^n$ in the case $1 \leq
\alpha <2$, some hypothesis are required. For two probability measures $\mu$,
$\nu$ in $\Omega$, we define the entropy $H(\mu|\nu)$ of $\mu$  with respect to
$\nu$ by $H(\mu|\nu)=\int d\mu/d\nu \log (d\mu/d\nu) d\nu$ if $\mu$ is
absolutely continuous with respect to $\nu$ and $H(\mu|\nu) =\infty$ if $\mu$ is
not absolutely continuous with respect to $\nu$.  

Consider now a sequence of measures $\{\mu^n\}_n$ associated to $u_0$ and such
that $\mu^n\{\eta(0)=1\}=1$. For each $n$ fixed, we tag the particle initially
at the origin, and we denote now by $X_t^n$ its position at time $tn^\alpha$.  
We have an invariance principle for $X_t^n$:

\begin{theorem}
\label{t4}
Assume that $0<\alpha <1$. Under the previous conditions on the initial measures
$\{\mu^n\}_n$, the process $X_t^n$ converges in distribution to the process of
independent increments $\mc Z_t$, characterized by the martingale problem 
\begin{equation}
\label{ec2}
M^\beta_t = \exp \big\{ i \beta \mc Z_t + \int_0^t \int q(x)\big( 1- e^{i\beta
x} \big)\big(1-u(s,x+\mc Z_s)\big) dx ds\big\} 
\end{equation}
is a martingale for any $\beta \in \bb R^d$, where $u(t,x)$ is the solution of
the hydrodynamic equation (\ref{echid}). 
Fix a density $\rho \in (0,1)$. For $1 \leq \alpha <2$, the same result is true
under the additional hypothesis 
\begin{equation}
\label{eb}
\sup_n n^{-d} H(\mu^n|\mu_\rho) < +\infty.
\end{equation} 
\end{theorem}

The entropy bound (\ref{eb}) plus the fact that $\{\mu^n\}_n$ is associated to the profile $u_0$ imply that $u_0-\rho$ is absolutely integrable. In fact, for a product measure $\mu^n$ such that $\mu^n\{\eta(x)=1\}=u_0(x/n)$, the entropy bound \ref{eb} is satisfied if and only if $\int|u_0(x)-\rho|dx$ is finite. A large deviations argument shows that for a given profile $u_0$, the entropy of $\{\mu^n\}_n$ with respect to $\mu_\rho$ is minimized by a sequence of measures of product form.

The process $\mc Z_t$ turns out to be a Markov process. The evolution of $\mc Z_t$ can be understood as follows. Put initially a particle at the origin. Then, this particle tries to move following the jumps of a L\'evy process $Z_t$. Each jump is accomplished with probability $1-u(t,x)$, where $x$ is the point of arrival. This is done independently for each jump. In particular, for $u_0=\rho$ we recover Theorem \ref{t1}.

\section{Invariance principle for $X_t^n$: the equilibrium case}
\label{s5}
For each $z \in \bb Z^d_*$ and each $\theta \in \bb R$, the process 
\[
 M_t^z = \exp\{i \theta N_t^z + (1-e^{i\theta}) \int_0^t
 p(z)\big(1-\xi_s(z)\big)ds\} 
\]
is a mean-one complex martingale. Since the jumps of $M_t^z$, $M_t^{z'}$ for $z
\neq z'$ are all different with probability one, for any $\beta \in \bb R^d$,
choosing $\theta = \beta z /n$, 
\begin{equation}
\label{ec4}
 M_t^{\beta,n} = \exp\big\{i\beta X_t^n + \sum_{z \neq 0} p(z) (1-e^{i\beta
 z/n})\int_0^{tn^\alpha}p(z)\big(1-\xi_s(z)\big)ds\big\} 
\end{equation}
is an exponential martingale with $E[M_t^{\beta,n}]=1$ for any $t, \beta, n$. We
follow the standard approach to prove Theorem \ref{t1}. First we prove
convergence of the finite-dimensional distributions of $X_t^n$ to the
corresponding distributions of $(1-\rho)Z_t$. We do this in Section \ref{s1}.
Then in Section \ref{s2} we prove tightness of the distributions of
$\{X_t^n\}_n$ in the Skorohod space $\mc D([0,T],\bb R^d)$, from which
convergence follows.  

\subsection{Convergence for $t>0$ fixed}
\label{s1}
In this section we prove the following

\begin{theorem}
\label{t2}
For each $t >0$, $X_t^n$ converges in distribution to $(1-\rho)Z_t$. 
\end{theorem}
\begin{proof}
It is enough to prove pointwise convergence of the corresponding characteristic
functions. In other words, we just need to prove that 
\[
\lim_{n \to \infty} E[\exp\{i\beta X_t^n\}] = E[\exp\{i\beta (1-\rho)Z_t\}] =
\exp\{-(1-\rho)t \psi(\beta)\} 
\]
for any $\beta \in \bb R^d$. Let us define 
\[
v_n(\xi) = n^\alpha \sum_{z \in \bb Z^d_*} p(z) (1-e^{i\beta
z/n})\big(1-\xi(z)\big)  
\]
and notice that 
\[
\int v_n d \nu_\rho = (1-\rho) n^\alpha \sum_{z \in \bb Z^d_*} p(z) (1-e^{i\beta
z/n}) = \frac{1-\rho}{n^d} \sum_{z \in \bb Z^d_*} q(z/n) (1-e^{i\beta z/n}). 
\]

This last sum is a Riemann sum for the integral $(1-\rho)\int (1-e^{i\beta
u})q(u) du$, and therefore 
\[
\lim_{n \to \infty} \int v_n d\nu_\rho = (1-\rho)\psi(\beta).
\]

The function $\Re(1-e^{i\beta u})=1-\cos(\beta u)$ is of quadratic order around
0. In particular, the integral $\int \Re(1-e^{i\beta u}) q(u) du$ is absolutely
summable, and $\sup_n||\Re(v_n)||_\infty < +\infty$. 
Since $|e^{i\omega}|=1$ for any $\omega \in \bb R$, we see that
$\exp\{\int_0^{tn^\alpha} n^{-\alpha}v_n(\xi_s)ds\}$ is uniformly bounded in
$n$. The martingale $M_t^{\beta,n}$ can be written in terms of $v_n$ and
$X_t^n$: 
\[
M_t^{\beta,n} = \exp\Big\{i\beta X_t^n + \frac{1}{n^\alpha}\int_0^{tn^\alpha}
v_n(\xi_s) ds\Big\}. 
\]

As we will see, Theorem \ref{t1} is a simple consequence of the following Lemma:
\begin{lemma}
\label{l1}
For any $t > 0$,
\[
\lim_{n \to \infty} \frac{1}{n^d} \int_0^{t n^\alpha} v_n(\xi_s) ds = (1-\rho) t
\psi(\beta) \text{ in probability.}
\]
\end{lemma}

In fact, since $E[M_t^{\beta,n}]=1$ we have
\[
\big|E[e^{i\beta X_t^n} - e^{-(1-\rho)t\psi(\beta)}]\big| 
	\leq E\big[\big|1-e^{n^{-\alpha} \int_0^{t n^\alpha} v_n(\xi_s)ds -
	(1-\rho) t \psi(\beta)}\big|\big],
\]
and this last quantity goes to zero by Lemma \ref{l1} plus $\sup_n||\Re
(v_n)||_\infty < +\infty$.
\end{proof}

Now we turn to the proof of the Lemma.
\begin{proof}[Proof of Lemma \ref{l1}]
Define $\bar v_n = v_n - \int v_n d\nu_\rho$ and $u_n = \Re(\bar v_n)$, $w_n =
\Im(\bar v_n)$. Notice that $w_n$ has a worse singularity at the origin than
$u_n$, and both functions have the same decay at infinity. Therefore, we will
prove the lemma just for $w_n$, the proof for $u_n$ being easier. The variance
of $w_n$ can be explicitly computed to obtain the bound
\[
\int (w_n)^2 d\nu_\rho \leq \frac{C(q,\rho)}{n^d} \Big\{ 1+\beta^2
\int_{1/n}^1\frac{du}{u^{d+2\alpha-1}}\Big\}.
\] 

In particular, $\int (w_n)^2 d\nu_\rho \leq C(q,\rho)n^{2\alpha-2}$ and for $0<
\alpha <1$, $\int (w_n)^2d\nu_\rho \to 0$ as $n \to \infty$, and the lemma is
trivial (there is no need of time integration). When $1\leq \alpha <2$ an extra
argument is needed. The idea is the following. For sites $x$ near the origin,
jumps of opposite sign cancel each other, as in the case of a diffusive system.
Since the transition rate $p(\cdot)$ has infinite second moment, these
cancellations do not hold for large $x$.
Introduce a truncation around 0 of size $n^\gamma$. If $\gamma$ is big
enough, by the previous argument the variance of the truncated function vanishes
as $n \to \infty$. In the other hand, if $\gamma$
is small enough, the central limit variance of small jumps vanishes too. We will
prove that $\gamma$ can be chosen in such a way that both things happen
simultaneously.

For $l \in \bb N$, define 
\[
w_n^l(\xi) =  n^\alpha \sum_{|z| \leq l} p(z) \sin(\beta
z/n)\big(\rho-\xi(z)\big).
\]

We have the bound
\[
\int (w_n - w_n^l)^2 d\nu_\rho \leq \frac{C(q,\rho)}{n^d} \Big\{ 1+\beta^2
\int_{l/n}^1\frac{du}{u^{d+2\alpha-1}}\Big\},
\]
and taking $l = n^\gamma$, $\int (w_n -w_n^{n^\gamma})d\nu_\rho \to 0$ as $n \to
\infty$ for $\gamma > (2\alpha -2)/(d+2\alpha-2)$. 
Of course, in order to have
cancellations of small jumps, we need the time integral. We have the following
estimate, valid for any function in $\Omega_*$ \cite{KV}:
\begin{equation}
\label{ec1}
E\Big[ \sup_{t \leq T} \Big( \frac{1}{n^\alpha} \int_0^{tn^\alpha} w_n^l(\xi_s)
ds \Big)^2\Big] 
\leq \frac{20 T}{n^\alpha} ||w_n^l||^2_{-1}, \text{ where }
\end{equation}
\[
||w_n^l||^2_{-1} =: \sup_{ f \in \mc L^2(\nu_\rho)} \Big\{ 2\< f, w_n^l\>_\rho -
\< f, -L f\>_\rho\Big\}
\]
and $\<\cdot,\cdot\>_\rho$ denotes the inner product in $\mc L^2(\nu_\rho)$.
Since the time $t$ is fixed in this section, we do not need the supremum inside
the expectation, but it will be needed later. In order to obtain an estimate of  
$2\< f, w_n^l\>_\rho$ in terms of $\<f,-Lf\>_\rho$, we first do some
manipulations:
\begin{equation}
\label{ec5}
\begin{split}
2\< f, w_n^l\> &= n^\alpha \sum_{|z| \leq l} p(z) \sin(\beta z/n) \<\xi(-z)
-\xi(z),f\>_\rho\\
	&=n^\alpha \sum_{|z| \leq l} p(z) \sin(\beta z/n)
		\<\xi(z), f(\xi^{z,-z})-f(\xi)\>_\rho \\ 
	&\leq n^\alpha \sum_{|z| \leq l} p(z) \sin(\beta z/n)
		\Big\{ \frac{\lambda(z)}{2} \<\xi(z),\xi(z)\>_\rho +\\
		&\qquad +\frac{1}{2\lambda(z)} \int
		\big(f(\xi^{z,-z})-f(\xi)\big)^2 d\nu_\rho\Big\}, \\
\end{split}
\end{equation}

where we have used the symmetry of $p(\cdot)$ in the first line, the invariance of $\nu_\rho$ under 
$\xi \to \xi^{z,-z}$ in the second line and Cauchy-Schwartz inequality in the third line.
After a change of variables, the inner product $\<f,-Lf\>_\rho$ can be written
as
\[
\<f,-Lf\>_\rho = \frac{1}{4} \sum_{x, y \in \bb Z^d_*} p(y-x) \int
\Big(f(\xi^{x,y})-f(\xi)\Big)^2 d\nu_\rho.
\]

Therefore, a simple way to estimate $||w_n^l||_{-1}^2$ is choose $\lambda(z)$ in
such a way that the sum bounding $\<w_n^l,f\>_\rho$ can be balanced by $\< f,
-Lf\>_\rho$. In other words, we choose
\[
n^\alpha p(z) \sin(\beta z/n) \frac{1}{2\lambda(z)} = \frac{1}{4} p(2z) =
\frac{1}{2^{d+\alpha+2}} p(z),
\]
with $\lambda(z) = 2^{d+\alpha+1}n^\alpha \sin(\beta z/n)$. We obtain in this
way the bound 
\[
||w_n^l||_{-1}^2 \leq 2^{d+\alpha} n^\alpha \frac{1}{n^d} \sum_{|z| \leq l}
q(z/n) \sin^2(\beta z/n).
\]

The sum is a Riemann sum of the integral $\int_0^{l/n} q(u) \sin^2(\beta u) du$.
For $l/n = o(1)$, this integral is of order $(l/n)^{2-\alpha}$. Putting this
estimate into (\ref{ec1}), we obtain the bound 
\[
E\Big[ \sup_{t \leq T} \Big( \frac{1}{n^\alpha} \int_0^{tn^\alpha} w_n^l(\xi_s)
ds \Big)^2\big] \leq c t(l/n)^{2-\alpha}.
\]
for some constant $c=c(q,d)$. Choosing $l=n^\gamma$, we see that
$(l/n)^{2-\alpha}$ goes to 0 as $n \to \infty$  
for any $\gamma <1$.
Therefore, taking $\gamma$ between $1- d/(d+2\alpha-2)$ and $1$, we have proved
that $E[(n^{-\alpha} \int_0^{tn^\alpha} w_n(\xi_s) ds)^2] \to 0$ as $n \to
\infty$. 
\end{proof}

\begin{remark}
Considering the martingales $M_{t,s}^{\beta,n} = M_{t+s}^{\beta,n}
(M_t^{\beta,n})^{-1}$, it is easy to prove that the increments of $X_t^n$ are
conditionally independent and identically distributed. It follows that for any
finite sequence of times $t_1<...<t_k$ we have
\[
(X_{t_1}^n,...,X_{t_k}^n) \to (1-\rho)(Z_{t_1},...,Z_{t_k}) \text{
in distribution.}
\]
\end{remark}

\subsection{Tightness of the sequence $\{X_t^n\}_n$}
\label{s2}
Fix some positive time $T$ and let us denote by $Q^n$ the distribution of
$X_t^n$ in the Skorohod space $\mc D([0,T],\bb R^d)$ of c\`adl\`ag trajectories.
Let us denote by $\mc T_T$ the set of stopping times bounded by $T$, with respect to the
filtration $\mc F_t$ associated to $\xi_t$. We prove
tightness of the sequence $\{Q^n\}_n$ by using Aldous' criterion:

\begin{proposition}[Aldous]

\end{proposition}
Let $\{Q^n\}_n$ be a family of probability distributions in $\mc D([0,T], \bb
R^d)$. The family $\{Q^n\}_n$ is relatively compact on the set of probability
distributions in $\mc D([0,T], \bb R^d)$ if
\begin{itemize}
\item[i)] The family $\{Q^n(\omega(0) \in \cdot)\}_n$ of distributions in $\bb R^d$
is relatively compact.
\item[ii)] For any $\epsilon >0$,
\[
\lim_{\delta \to 0} \limsup_{n \to \infty} \sup_{\substack{\tau \in \mc
T_T\\\gamma \leq \delta}} Q^n(|\omega(\tau+\gamma) -\omega(\tau)|>\epsilon) =0.
\]
\end{itemize}

Since $X_0^n \equiv 0$, the first condition is trivially satisfied. To prove
$ii)$, we split $X_t^n$ in two pieces: let $l=l(n)$ be a sequence to be chosen
later and define
\[
X_t^{n,1} = \frac{1}{n}\sum_{|z| \leq l} z N_{tn^\alpha}^z, 
X_t^{n,2} = \frac{1}{n}\sum_{|z| > l} z N_{tn^\alpha}^z.
\]

In other words, $X_t^{n,1}$ corresponds to the small jumps of $X_t^n$ and
$X_t^{n,2}$ corresponds to the long jumps of $X_t^n$. Tightness of $X_t^n$ will
follow from tightness of each one of these processes. Let us start with
$X_t^{n,2}$. We have the estimate
\[
P(|X_{\tau+\gamma}^{n,2} - X_\tau^{n,2}| > \epsilon) \leq \frac{1}{2\epsilon}
\int_{|\beta| \leq 1/\epsilon} \big| 1- E[e^{i\beta(X_{\tau+\gamma}^{n,2} -
X_\tau^{n,2})}]\big|d\beta.
\]

Define $M_t^{n,2}$ in the natural way. The process $M_{\tau+t}^{n,2}
(M_{\tau}^{n,2})^{-1}$ is a $\mc F_{\tau+t}$-martingale, and therefore
\[
E\big[ \exp\big\{ i\beta(X_{\tau+t}^{n,2} - X_\tau^{n,2}) +n^{-\alpha}
\int_{\tau n^\alpha}^{(\tau+t)n^\alpha}
\big(v_n(\xi_s)-v_n^l(\xi_s)\big)ds\big\}\big] =1,
\]
where $v_n^l$ is the truncation of $v_n$, defined as we did for $w_n^l$.
In particular,
\[
\big| 1- E[e^{i\beta(X_{\tau+\gamma}^{n,2} -
X_\tau^{n,2})}]\big| \leq E \big| 1 - \exp\big\{ n^{-\alpha}
\int_{\tau n^\alpha}^{(\tau+\gamma)n^\alpha}
\big(v_n(\xi_s)-v_n^l(\xi)\big)ds\big\} \big|.
\]

Choosing $l=n$, we see that $v_n -v_n^l$ is uniformly bounded in $\xi$ and uniformly bounded in
compact sets in $\beta$. Therefore, the previous integral involving $v_n -v_n^n$
can be bounded by $C(\epsilon)\delta$ for $\beta \in [-1/\epsilon, 1/\epsilon]$ uniformly
in $n$, and condition $ii)$ is satisfied by $X_t^{n,2}$. Notice that in the case $0<\alpha <1$, there is no need of truncation: $v_n$ is already bounded. In particular the proof of tightness for $0<\alpha<1$ is complete. We assume then $1\leq \alpha <2$.

As pointed out before, the exponential martingales $M_t^{\beta,n}$ were
introduced since $X_t^n$ has no bounded second moments. This is not longer a
problem for $X_t^{n,1}$. For each $z \in \bb Z^d_*$, the process $M_t^z = N_t^z
- \int_0^t p(z)\big(1-\xi_s(z)\big) ds$ is a martingale of quadratic variation
$\int_0^t p(z)\big(1-\xi_s(z)\big) ds$. As before, the martingales $\{M_t^z,z
\in \bb Z_*^d\}$ are mutually orthogonal. We can write $X_t^{n,1}$ as
\[
X_t^{n,1} = \frac{1}{n} \sum_{|z| \leq n} z N_{tn^\alpha}^z 
	= \frac{1}{n} \sum_{|z| \leq n} z M_{tn^\alpha}^z
		+ \frac{1}{n}\sum_{|z| \leq n} z p(z) \int_0^{t n^\alpha}
		\big(1-\xi_s(z)\big) ds.
\] 

Therefore, $X_t^{n,1} = M_t^{n,1} + n^{-\alpha} \int_0^{tn^\alpha}
v_{n,1}(\xi_s) ds$, where $M_t^{n,1}$ is a martingale of quadratic variation
\[
\<M_t^{n,1} \> = \sum_{|z| \leq n} \frac{z^2}{n^2} p(z) \int_0^{tn^\alpha}
v_{n,1}^n(\xi_s) ds, 
\]
and we have defined 
\[
v_{n,1}^l(\xi) = \frac{1}{2n^d} \sum_{|z| \leq l} \frac{z}{n}
q(z/n)\big(\xi(-z)-\xi(z)\big).
\]

In this last identity, we have used the symmetry of $q(\cdot)$. Notice that
$\<M_t^{n,1}\> \leq tn^{-d} \sum_{z} (z/n)^2 q(z/n)$. By Aldous' criterion and
the optional sampling theorem, the sequence of processes $\{M_t^{n,1}\}_n$ is
tight:
\begin{align*}
P\big(\big|M_{\tau+\gamma}^{n,1}-M_{\tau}^{n,1}\big|>\epsilon\big) &\leq
\frac{1}{\epsilon^2}
E\big[\big(M_{\tau+\gamma}^{n,1}-M_{\tau}^{n,1}\big)^2\big]\\ 
	&\leq \frac{1}{\epsilon^2}
	E\big[\<M_{\tau+\gamma}^{n,1}\>-\<M_{\tau}^{n,1}\>\big]\\
	& \leq \frac{1}{\epsilon^2} \frac{\delta}{n^d} \sum_{|z/n| \leq 1}
	(z/n)^2 q(z/n).
\end{align*}

Since $x^2 q(x)$ is absolutely integrable in the unit ball, we conclude that
\[
\lim_{\delta \to 0} \limsup_{n \to \infty} \sup_{\substack{\tau \in \mc
T_T\\\gamma \leq \delta}}
P\big(\big|M_{\tau+\gamma}^{n,1}-M_{\tau}^{n,1}\big|>\epsilon\big) =0. 
\]

We are left to prove tightness for the process $n^{-\alpha} \int_0^{tn^\alpha} v_{n,1}^n(\xi_s) ds$. 
Remember the computations after Lemma \ref{l1}. Notice the similarity between the functions 
$v_{n,1}$ and $w_n$. 
For $l = n^\theta$ with $\theta$ between $(2\alpha-2)/(d+2\alpha-2)$ and 1, 
$E[(v_{n,1}^n - v_{n,1}^l)^2]$ goes to 0 as $n \to \infty$. In the other hand,
\begin{align*}
E\Big[\Big(n^{-\alpha}\int_{\tau n^\alpha}^{(\tau+\gamma)n^\alpha}
	&\big(v_{n,1}^n(\xi_s)-v_{n,1}^l(\xi_s)\big)ds\Big)^2\Big] \leq \\
	&\leq \gamma E\Big[ n^{-\alpha}\int_{\tau n^\alpha}^{(\tau+\gamma)n^\alpha}
		\big(v_{n,1}^n(\xi_s)-v_{n,1}^l(\xi_s)\big)^2 ds\Big] \\
	&\leq \gamma T E\big[\big(v_{n,1}^n-v_{n,1}^l\big)^2\big] \xrightarrow{n \to \infty} 0.
\end{align*}

Following the same computations in the proof of Lemma \ref{l1}, we can obtain 
\[
\lim_{n \to \infty} E\Big[\sup_{t \leq T}\Big(\int_0^{t n^\alpha} v_{n,1}^l(\xi_s) ds\Big)^2\Big] =0,
\]
and in particular $n^{-\alpha} \int_0^{t n^\alpha} v_{n,1}^n(\xi_s)ds $ is tight, which end the proof of tightness for $\{Q^n\}_n$.

\section{Invariance principle for $X_t^n$: the non-equilibrium case}
\label{s6}
In this section we prove Theorem \ref{t4}. We follow the approach introduced in \cite{JLS}. 
That is, first we prove tightness in the Skorohod space for the sequence $\{X_t^n\}$. Therefore, the joint process $(\pi_t^n, X_t^n)$ is tight as well. In particular, the empirical density as seen by the tagged particle $\bar \pi_t^n$ is also tight, where
\[
\bar \pi_t^n(dx) =\frac{1}{n^d} \sum_{z \in \bb Z^d_*} \xi_t^n(z) \delta_{z/n}(dx).
\]

Take a subsequence for which each one of the previous processes converges.
Then, using the martingale representation (\ref{ec4}) of $X_t^n$ and the convergence of $\pi_t^n$, $\bar \pi_t^n$, $X_t^n$, we obtain a martingale characterization of the limit points of $X_t^n$. In fact, it turns out that $v_n(\xi_s)$ is a function of the empirical measure. Notice the similarity between (\ref{ec4}) and (\ref{ec2}). We finish the proof by an uniqueness result for the process satisfying (\ref{ec2}). 

This program can be accomplished for $0 <\alpha <1$ in a simple way. 
As we saw in Section \ref{s2}, when $1\leq \alpha <2$ small jumps cancel each other in a non-trivial way. The needed estimates are more intricate out of equilibrium. In the other hand, since $q(x)(1-e^{i\beta x})$ is not uniformly integrable, $v_n(\xi_s)$ is not longer a function of the empirical measure. A cut-off argument to take into account cancellations of small jumps is needed again.

\subsection{Tightness and characterization of limit points: $0 < \alpha <1$}
\label{s3}
We start proving tightness of $\{X_t^n\}$. Recall \ref{ec4}:
\[
M_t^{\beta,n} = \exp\big\{ i \beta X_t^n + \int_0^{t} v_n(\xi_s^n) ds\big\} 
\]
is a martingale. Therefore, for any stopping time $\tau$ bounded by $T$ and any $\gamma>0$,
\begin{align*}
P(|X_{\tau+\gamma}^{n} - X_\tau^{n}| > \epsilon) 
	&\leq \frac{1}{2\epsilon} \int_{|\beta| \leq 1/\epsilon} 
	\big| 1- E[e^{i\beta(X_{\tau+\gamma}^{n} - X_\tau^{n})}]\big|d\beta \\
	& \leq \frac{1}{2\epsilon} \int_{|\beta| \leq 1/\epsilon} 
	\Big| 1- E \exp\Big\{\int_\tau^{\tau+\gamma} v_n(\xi_s) ds\Big\}\Big| d\beta.
\end{align*}

Since $v_n$ is uniformly bounded in $n$, $\xi$ and for $\beta$ in compact sets, we conclude that $X_t^n$ is tight.
In general, tightness of the vector $(X_t^n,\pi_t^n)$ can not be obtained from the tightness of each component. However, since $\pi_t^n$ converges in probability to a deterministic limit, we have that $(X_t^n,\pi_t^n)$ is tight in $\mc D([0,\infty), \bb R^d \times \mc M_+(\bb R^d))$. Take a limit point $\mc Z_t$ of $X_t^n$ and take a subsequence (still denoted by) $n$ such that $(X_t^n, \pi_t^n(dx))$ converges to $(\mc Z_t, u(t,x) dx)$. Then, for any continuous, bounded function $G: \bb R^d \to \bb R$ we have
\[
\lim_{n \to \infty} \int _0^t \int G(x) \bar \pi_s^n(dx) ds= \int_0^t \int_{\bb R^d} G(x) u(s,x+\mc Z_s) dx ds
\]
in distribution. Notice that $v_n(\xi_s^n) = \int q(x)(1-e^{i\beta x})\big(1- \bar \pi_s^n\big)(dx)$, where $1- \bar \pi_s^n$ is the empirical measure defined with $1-\xi_t^n(x)$ instead of $\xi_t^n(x)$. The function $q(x) (1-e^{i \beta x})$ is not continuous at $x=0$, but it is absolutely integrable. In particular, there exists a sequence $\{G_n\}_n$ of continuous functions with bounded support such that $G_n(x)$ converges to $q(x)(1-e^{i \beta x})$ in $\mc L^1(\bb R^d)$. Moreover, we can take $G_n = q(x)(1-e^{i\beta x})$ in a ring of inner radius $1/n$ and outer radius $n$, centered at the origin.  Since the number of particles per site is bounded, we conclude that
\[
\lim_{n \to \infty} \int_0^t v_n(\xi_s^n) ds = \int_0^t \int_{\bb R^d} q(x) (1-e^{i\beta x})\big(1-u(s,x+\mc Z_t)\big)dx ds.
\]

Since everything in (\ref{ec4}) is bounded, by the Dominated Convergence Theorem we conclude that 
\[
M_t^\beta = \exp\big\{ i\beta \mc Z_t + \int_0^t \int q(x) (1-e^{i\beta x}) \big(1-u(s,x+\mc Z_s)\big)dx ds\big\}
\]
 is a martingale for any $\beta \in \bb R^d$.

\subsection{Tightness and characterization of limit points: $1\leq \alpha <2$}

As we have learned from Section \ref{s2}, when $1 \leq \alpha <2$, an argument to take into account cancellation of small jumps is needed. First we observe that the arguments of the previous section allow us to prove tightness of the truncated process
\[
X_t^{n,\epsilon,2} = \frac{1}{n} \sum_{|z| \geq \epsilon n} z N_{tn^\alpha}^z
\]
for each $\epsilon >0$.
Therefore, we are left with the problem of tightness for $X_t^{n,\epsilon}= X_t^n - X_t^{n,\epsilon,2}$. The jumps of $X_t^{n,\epsilon}$ are bounded by $\epsilon$. Remember the decomposition
\[
X_t^{n,\epsilon} = M_t^{n,\epsilon} + \int_0^t v_{n,1}^{\epsilon n}(\xi_s^n) ds.
\]

The proof of tightness for $M_t^{n,\epsilon}$ follows as in the equilibrium
case. Therefore, we are left with the problem of tightness for the integral part
$\int_0^t v_{n,1}^{\epsilon n}(\xi_s^n)ds$. The whole point here is to estimate
integrals  of the form $(z/n)^{-1}\int_0^t \{\xi_s^n(z)-\xi_s^n(-z)\}ds$. The same
problem arises in the hydrodynamic limit of non-gradient systems \cite{VY},
\cite{Qua}. We start with an exponential bound on the norm of the integral term:

\begin{lemma}
\label{l2}
There exist constant $c>0$ and function $\lambda(\epsilon)$ with $\lambda(\epsilon) \to 0$ as $\epsilon \to 0$, such that
\[
\bb E_\rho \Big[ \exp\Big\{ n^d \Big| \int_{t_0}^{t_1} v_{n,1}^{\epsilon n}(\xi_s^n) ds\Big|\Big\}\Big] 
	\leq c \exp\big\{n^d(t_1-t_0)\lambda(\epsilon)\big\},
\]
where $\bb E_\rho$ corresponds to the expectation with respect to the process $\xi_t^n$ starting from $\mu_\rho$.
\end{lemma}

\begin{proof}
We can get rid of the modulus by observing that $e^{|x|} \leq e^x+e^{-x}$. By Feynman-Kac formula and the variational expression for the largest eigenvalue of a symmetric operator, we have the estimate
\[
\log \bb E_\rho \Big[ \exp\Big\{ n^d \int_{t_0}^{t_1} v_{n,1}^{\epsilon n}(\xi_s^n) ds\Big\}\Big]
	\leq (t_1-t_0)n^d \sup_f \big\{ \< v_{n,1}^{\epsilon n}, f^2\>_\rho - \<f, -L f\>_\rho\big\},
\]
where the supremum is over functions $f: \Omega_* \to \bb R^d$ with $\<f,f\>_\rho=1$. Computations very similar to \ref{ec5} allow us prove that the right-hand side of the previous expression is bounded by
\[
(t_1-t_0) c\sum_{|z|\leq \epsilon n} (z/n)^2 q(z/n) \leq n^d(t_1-t_0) \lambda(\epsilon),
\]
where we have used the fact that $x^2 q(x)$ is absolutely integrable in the unit ball. 
\end{proof}

The following result is a restatement of (\cite{KL}, Corollary 7.6.6.4), and it
is a (not so immediate) consequence of Lemma \ref{l2}. 

\begin{proposition}
For each $T, \epsilon>0$ there exists a constant $c(\epsilon, T)$ with
$c(\epsilon,T) \to 0$ as $\epsilon \to 0$ such that 
\[
\sup_{n >0} \bb E^n\Big[ \sup_{\substack{|t-s| < \delta\\s,t \leq T}}
\Big|\int_s^t v_{n,1}^{\epsilon n}(\xi_s^n) ds \Big| \Big] \leq c(\epsilon,T)
\sqrt \delta \log(1/\delta). 
\]
\end{proposition}

\begin{remark}
It is in the proof of the previous proposition that the entropy bound (\ref{eb})
is needed. The proof of Theorem \ref{t4} is otherwise independent of (\ref{eb}).
See (\cite{KL}, Section 7.6) for more details. 
\end{remark}

We therefore have proved tightness for $\int_0^t v_{n,1}^{\epsilon n}(\xi_t^n)
ds$, from where tightness of $X_t^n$ follows. Notice that we still did not make
use of the fact that $c(\epsilon,T) \to 0$.  
As in the previous section, tightness of $(\pi_t^n, X_t^n)$ follows. However,
now the function $q(x)(1-e^{i\beta x})$ is not absolutely integrable. We
introduce a truncation in order to overcome this problem. Define
$Z_t^{n,\epsilon} = X_t^n -X_t^{n,\epsilon}$. Tightness of $Z_t^{n,\epsilon}$
follows as before. Let $\mc Z_t^\epsilon$  
be a limit point of $Z_t^{n,\epsilon}$, and take a subsequence such that $X_n^t$
converges to some process $\mc Z_t$ and $Z_t^{n,\epsilon}$ converges to $\mc
Z_t^\epsilon$. Define $q_\epsilon(x) = q(x) {\bf 1}(|x| \geq \epsilon)$. Then,  
\[
M_t^\beta = \exp\big\{ i\beta \mc Z_t^\epsilon + \int_0^t \int q_\epsilon(x) (1-e^{i\beta x}) \big(1-u(s,x+\mc Z_s)\big)dx ds\big\}
\]
is a martingale for any $\beta \in \bb R^d$. It is clear that for each $n>0$,
$Z_t^{n,\epsilon}$ converges to $X_t^n$ as $\epsilon \to 0$. By Lemma \ref{l2},
this convergence is uniform in $n$. Here we need the fact that $c(\epsilon,T) \to 0$. Therefore, we can exchange the limits $n \to \infty$ and $\epsilon \to 0$ to prove that $\mc Z_t^\epsilon$ converges to $\mc Z_t$ as $\epsilon \to 0$. By the same arguments, we conclude that
\[
M_t^\beta = \exp\big\{ i\beta \mc Z_t+ \int_0^t \int q(x) (1-e^{i\beta x}) \big(1-u(s,x+\mc Z_s)\big)dx ds\big\}
\]
is a martingale for any $\beta \in \bb R^d$.

\subsection{Uniqueness of the martingale problem}

In the previous section we proved that all the limit points of $X_t^n$ satisfy the martingale problem (\ref{ec2}). Denote by $\mc Z_t$ any of such limit points. It follows that (see 4.2.1 and 1.2.8 in \cite{SV}) $\mc Z_t$ satisfies the martingale problem
\[
F(t, \mc  Z_t) - \int_0^t \big\{\partial_t +\mc L_s \big\} F(s,\mc Z_s) ds
\]
is a martingale for any bounded function $F \in \mc C^{1,2}([0,\infty)\times \bb R^d)$, where $\mc L_t$ is the integral operator
\begin{align*}
\mc L_t F(x) &=  \int_{\bb R^d} q(y)\{F(x+y)-F(x)\}\{1-u(t,x+y)\}dy,
\end{align*}
where we have to consider the principal value of the integral. 
This operator corresponds to the generator of a time-inhomogeneous process of independent increments. Uniqueness for this martingale problem can be obtained by adapting the arguments of \cite{SV} to this setting (see also \cite{EK}). 

\section*{acknowledgments}
The author would like to thank valuable discussions with L. Silvestre and S.R.S Varadhan, and the kind hospitality at the Courant Institute, where this work was initiated.

\def\cprime{$'$}


\begin{thebibliography}{10}

\bibitem{EK}
Stewart~N. Ethier and Thomas~G. Kurtz.
\newblock {\em Markov processes}.
\newblock Wiley Series in Probability and Mathematical Statistics: Probability
  and Mathematical Statistics. John Wiley \& Sons Inc., New York, 1986.
\newblock Characterization and convergence.

\bibitem{Fer}
P.~A. Ferrari.
\newblock Limit theorems for tagged particles.
\newblock {\em Markov Process. Related Fields}, 2(1):17--40, 1996.
\newblock Disordered systems and statistical physics: rigorous results
  (Budapest, 1995).

\bibitem{JL1}
M.~D. Jara and C.~Landim.
\newblock Nonequilibrium central limit theorem for a tagged particle in
  symmetric simple exclusion.
\newblock {\em Ann. Inst. H. Poincar\'e Probab. Statist.}, 42(5):567--577,
  2006.

\bibitem{JL2}
M.~D. Jara and C.~Landim.
\newblock Quenched nonequilibrium central limit theorem for a tagged particle
  in symmetric simple exclusion.
\newblock {\em To appear in Ann. Inst. H. Poincar\'e Probab. Statist.}, 2006.
\newblock Available online at http://fr.arxiv.org/abs/math/0603653.

\bibitem{JLS}
M.~D. Jara, C.~Landim, and S.~Sethuraman.
\newblock Nonequilibrium fluctuations for a tagged particle in mean-zero
  one-dimensional zero-range processes.
\newblock Preprint.
\newblock Available online at http://fr.arxiv.org/abs/math/0703226.

\bibitem{KV}
C.~Kipnis and S.~R.~S. Varadhan.
\newblock Central limit theorem for additive functionals of reversible {M}arkov
  processes and applications to simple exclusions.
\newblock {\em Comm. Math. Phys.}, 104(1):1--19, 1986.

\bibitem{KL}
Claude Kipnis and Claudio Landim.
\newblock {\em Scaling limits of interacting particle systems}, volume 320 of
  {\em Grundlehren der Mathematischen Wissenschaften [Fundamental Principles of
  Mathematical Sciences]}.
\newblock Springer-Verlag, Berlin, 1999.

\bibitem{LOV}
C.~Landim, S.~Olla, and S.~R.~S. Varadhan.
\newblock Asymptotic behavior of a tagged particle in simple exclusion
  processes.
\newblock {\em Bol. Soc. Brasil. Mat. (N.S.)}, 31(3):241--275, 2000.

\bibitem{Lig}
Thomas~M. Liggett.
\newblock {\em Interacting particle systems}, volume 276 of {\em Grundlehren
  der Mathematischen Wissenschaften [Fundamental Principles of Mathematical
  Sciences]}.
\newblock Springer-Verlag, New York, 1985.

\bibitem{Qua}
Jeremy Quastel.
\newblock Diffusion of color in the simple exclusion process.
\newblock {\em Comm. Pure Appl. Math.}, 45(6):623--679, 1992.

\bibitem{Rez}
Fraydoun Rezakhanlou.
\newblock Propagation of chaos for symmetric simple exclusions.
\newblock {\em Comm. Pure Appl. Math.}, 47(7):943--957, 1994.

\bibitem{SVY}
Sunder Sethuraman, S.~R.~S. Varadhan, and Horng-Tzer Yau.
\newblock Diffusive limit of a tagged particle in asymmetric simple exclusion
  processes.
\newblock {\em Comm. Pure Appl. Math.}, 53(8):972--1006, 2000.

\bibitem{SV}
Daniel~W. Stroock and S.~R.~Srinivasa Varadhan.
\newblock {\em Multidimensional diffusion processes}, volume 233 of {\em
  Grundlehren der Mathematischen Wissenschaften [Fundamental Principles of
  Mathematical Sciences]}.
\newblock Springer-Verlag, Berlin, 1979.

\bibitem{Var}
S.~R.~S. Varadhan.
\newblock Self-diffusion of a tagged particle in equilibrium for asymmetric
  mean zero random walk with simple exclusion.
\newblock {\em Ann. Inst. H. Poincar\'e Probab. Statist.}, 31(1):273--285,
  1995.

\bibitem{VY}
S.~R.~S. Varadhan and Horng-Tzer Yau.
\newblock Diffusive limit of lattice gas with mixing conditions.
\newblock {\em Asian J. Math.}, 1(4):623--678, 1997.

\end{thebibliography}
\end{document}